\documentclass[twoside,11pt]{amsart}

\title{Challenges to Some Philosophical Claims about Mathematics}
\author{Eliahu Levy}
\address{Department of Mathematics,
Technion -- Israel Institute of Technology,
Haifa 32000, Israel}
\email{eliahu@techunix.technion.ac.il}

\newcommand{\NOT}[1]{}
\newcommand{\pa}{\par\medskip}

\newcommand{\bN}{\mathbb{N}}
\newcommand{\bR}{\mathbb{R}}
\newcommand{\cP}{\mathcal{P}}

\newcommand{\Ch}{{\bf Challenge.\ }}

\begin{document}
\maketitle
\begin{abstract}
In this note some philosophical thoughts and observations about mathematics are expressed, arranged as challenges to some common claims.
\end{abstract}
\renewcommand{\thesection}{}

\section{}
For many of the ``claims'' and ideas in the ``challenges'' see the sources listed in the references.\pa

\subsection{Claim}
The Antinomies in Set Theory, such as the Russell Paradox, just show that people did not have a right concept about sets. Having the right concept, we get rid of any contradictions.\pa
\Ch
It seems that this cannot be honestly said, when often in ``axiomatic'' set theory the same reasoning that leads to the Antinomies (say to the Russell Paradox) is used to prove theorems -- one does not get to the contradiction, but halts before the ``catastrophe'' to get a theorem.  As if the reasoning that led to the Antinomies was not ``illegitimate'', a result of misunderstanding, but we really have a contradiction (antinomy) which we, somewhat artificially, ``cut'', by way of the axioms, to save our consistency.\pa

One may say that the phenomena described in the famous G\"odel's Incompleteness Theorem
are a reflection of the Antinomies and the resulting inevitability of an axiomatics not entirely parallel to intuition. Indeed, G\"odel's theorem forces us to be presented with a statement (say, the consistency of Arithmetics or of Set Theory) which we know we cannot prove, while intuition puts a ``proof'' on the tip of our tongue, so to speak (that's
how we ``know'' that the statement is true!), but which our axiomatics, forced
to deviate from intuition to be consistent, cannot recognize. (Note that, contrary to Roger Penrose, there seems to be no further mysterious ``human superior to computer'' issue here.)\pa

To further clarify the issue, let us consider an example of a ``high-level'' part of Set Theory -- allowing ourselves to consider the power set $\cP(A)$ of a set $A$, i.e.\ the set of all subsets of $A$\NOT{(a set $S$ is defined to be a subset of $A$ if all members of $S$ are also members of $A$)}. One's ``attitude'' to this may swing between at least two poles: 1.\ When we ``have'' $A$, we have all its members, so we have the subsets of $A$, these ``grouping some of the members'' (and these subsets may be identified with all conceivable properties that the members of $A$ may have), so there is no problem in the set $\cP(A)$ of all of them. 2.\ A set must be {\em given} to us, in particular as the set of all elements satisfying some property, and sets which happen to be subsets of $A$ (a thing which may
be very hard to check/verify) may ``come'' from anywhere, ``from far away'', and grouping them to a set $\cP(A)$ is not far from speaking about ``the set of all sets'' which features in the Antinomies.\pa

Note, however, that one may say that adopting such attitudes {\em should not matter for us when we are doing the mathematics}. However ``vague'' such ``high-level'' Set Theory may be to the taste of people like Hermann Weyl, one must concede that, as Cantor has shown us, this Set Theory {\em is} something mathematical -- it deals with unambiguous notions, in the sense of having an ``adequate language'' to speak about them (so that there will not be questions which we do not know how to answer not because of lack of knowledge but because the terms in the question are too ``vague''). Consequently, complicated logical constructions are safely meaningful so that logical reasoning can be set to work and show its power.\pa

Note, however, a caveat here. Such attitudes may lead one to altogether reject, say, allowing oneself to make this step of considering the power set, thus to be entirely prepared to see theorems and theories, as legitimate, that would lead to a contradiction if considered together with such step. In other words, to treat making such a step in the same way one is used to do with the considerations leading to the Antinomies, hence banned by way of the choice of axioms.
\pa

\subsection{Claim}
The natural numbers $\bN$ are something concrete, basic to the mental world of human beings.\pa
\Ch
Here I distinguish between the concept of a general natural number on the one hand, and our practice with ``small'' natural numbers on the other hand. The latter are, say, the number {\em five} of the fingers in a hand, but also bigger numbers that we say or write their name, made of an appropriate linguistic combination (of words or digits). We use them to count, and as a consequence have the notions of adding them, multiplying them etc. We are used to learn their properties from experience (e.g.\ the crucial fact that a set cannot be counted by both $n$ and $n+1$, or more generally by both $n$ and $n+n'$, $n'\ge1$).\pa

The set $\bN$ of all natural numbers, on the other hand, is an abstract theoretical
structure, consisting of abstract entities, which is introduced as an extension of the above ``small numbers'' so that the algebraic operations are everywhere defined. Its introduction is partly a matter of mathematical convenience, enhanced by the inconvenience resulting from the existence of many non-isomorphic discrete segments $\{1,2,\ldots,n\}$. Certainly, the elements of $\bN$ are not parts of our everyday experience either with sensible objects or with linguistic constructions. (Note that a child must be quite sophisticated to conceive the notion of number not ``given'' with a definite/specific name).\pa

Nevertheless, $\bN$ has a strong intuitive appeal. People are sure that they can imagine a discrete segment of $10^{50}$ elements -- the imaginative picture starts with some initial elements, a ``gap'' and some final ones, and is quite the same as that of $10^{100}$ elements. Similarly, people are sure they can imagine a sequence enumerated by $\bN$, much in the same way as they imagine $10^{50}$ elements, only here the picture starts with some initial elements followed by an ``infinite gap''.\pa

Some people say: $\bN$ is a basic ingredient of our mental world, just start with $1$ (or $0$) and repeat indefinitely the operation of adding $1$. But, this ``repeating indefinitely'' is ``repeating indexed by $\bN$''. If we are in a position to repeat indefinitely, this mean that we already have the (abstract, theoretical) $\bN$. Indeed, if we were ``given'' the discrete segment $J:=\{1,2,\ldots10^{10^{1000}}\}$ we could mean by
``indefinitely'' being indexed by $J$.%
\footnote{If I knew I would live $10^{10^{1000}}$ days, could I not ignore that there will be a last day?}\pa

Now, one postulates about $\bN$ properties similar to those with which one is accustomed
in the above ``small'' numbers. In particular, we want the validity of mathematical induction -- proof by induction and definition by induction. This requirement is, in a sense, an implicit definition of the (again, abstract and theoretical) $\bN$.\pa

\NOT{
The conviction of people in the validity of IV) is so strong that very often in the teaching of mathematics something is defined or proved for some first small numbers, and it is agreed that this defines or proves it for every $n$ in $\bN$, especially if a ``recursive'' way to construct such definition or proof for general $n$ is apparent. This is a standard procedure in Euclid's ``Elements''. It hides the role of the Axiom of Induction (or an equivalent one).
}

\subsection{Claim}
It is OK to define a natural number as a thing common to equipotent finite sets (two sets being called {\em equipotent} if there is a one-one correspondence between them).\pa
\Ch
For that to be OK, one should previously define the notion of ``finite sets'' in a manner independent of natural numbers. Needless to say, this is almost never done. Those who employ such an approach seem to try to make the impression that ``finite sets are any sets you may stumble on'', which is contradicted the moment the set $\bN$ is conceived.%
\footnote{The fact that this inconsistency is almost always unnoticed by teachers and pupils shows again how remote the notion of ``general'' natural numbers and of $\bN$ is, in fact, from the common concepts of, say, a child.}

\subsection{Claim}
People do experiments with natural numbers, to test conjectures like the Riemann Hypothesis. If these do not give 100\% certainty, at least they give very high certainty like in the natural sciences. Similarly, one can establish formulas by showing that they are correct to, say, thousands or millions of digits.\pa
\Ch
Let us observe that the seemingly ``logically neat'' $\bN$  conceals a kind of paradox/dilemma in one's attitude towards it, so to speak. Indeed, for the purposes of almost everyone except the student of $\bN$ ``for its own sake'', the natural numbers (and consequently real numbers, with respect both to size and to accuracy) that ``count'' are
those that one is ``going to work with practically'', thus are all bounded by some $M$, usually much less than, say, $10^{1000000}$. On the other hand, for any $M\in\bN$, however big, the ``overwhelming majority'' of members of $\bN$ are greater than $M$. Thus sampling among the numbers which we can practically calculate with to test ``empirically'' some conjecture involving all natural numbers (say, something like the Riemann Hypothesis) must be considered a very poor and biased sampling, yet for those for whom a natural number means a number with which they may sometime ``work'', such sampling must seem quite
proper.\pa

A (highly heuristic) justification for taking natural numbers up to some $M$ as a proper sample is the following: One proceeds from the belief that the ``complexity'' of, say, the Riemann Hypothesis is somehow finite, hence everything about it must be revealed by the natural numbers up to some $M$, where persistence of the observed ``experimental'' phenomena when one takes higher and higher $M$ is taken as ``empirical evidence'' that the
true ``level of complexity'' has been reached.\pa

Indeed, only by invoking such a belief in ``not too high complexity'' can people claim that, say, a computation that showed that two real expressions involving functions, series etc.\ are equal in, say, the first one million decimal digits, should lead us to believe that they are equal as real numbers in the mathematical sense.\pa

Also, the idea that for most people ``the natural numbers'' are ``really'' the elements of some ``finite'' initial interval $\{n\,|\,n\le M\}$ (and that carries over, of course, to the real numbers, again concerning both size and accuracy) sheds a new and strange light on the issue of ``infinity''.\pa

\subsection{Claim}
Sets are altogether different things from numbers. Also, while in the first two thirds of the 20th century people based mathematics on sets, today categories and functors arise as a preferable alternative.\pa
\Ch
Firstly, the system of natural numbers, the time-honored basis for mathematics, may be naturally viewed as a version of the system of ``sets of finite rank'', which differs
from Set Theory only in the lack of the Axiom of Infinity, thus a ``ban'' on infinite sets. (And recall that Euclid's ``Elements'' speaks of numbers almost as the same as (finite) sets of ``units''.) Since, as we saw above, both numbers and sets are quite abstract, there seems to be no justification to admit mathematics based on numbers but not mathematics based on sets, where, to say it again, the difference between the two boils down to the presence of something like the Axiom of Infinity.
\pa

As for categories and functors, they are defined as some structures made from sets (or classes). Even if we choose to start everything from them, still we shall have the class (or set) of objects, the set of morphisms between two objects, etc. This seems to lie well inside basing mathematics on sets (over and above the natural numbers) -- these to be dealt with lavishly, making them members of further sets etc. That is the legacy of Georg Cantor and his followers, who showed convincingly that mathematics is indeed naturally based on sets in this sense.\pa

But of course, the flavor of Set Theory and basing ``everything'' on structures made of sets, \`{a} la Bourbaki, is rather different from the spirit of basing everything on categories and functors. Similarly, even if ultimately some branch of mathematics can be ``reduced'' to a part of, say, Set Theory, the flavor of it and of studying it would usually be quite different from that of Set Theory.

\subsection{Claim}
It is quite significant whether and how mathematical notions, such as numbers, exist.\pa
\Ch
This could really be claimed only if one is ready to view everything said about, say,
$\pi$, as gibberish if $\pi$ does not ``respectfully'' ``exist'' (as in the anecdote about Kronecker saying to Lindemann: ``it is very nice that you proved that $\pi$ is transcendental, but you are talking about something that does not exist!'').\pa

It seems that if we are allowed to talk about, say, the number $2$ or $\pi$ or the set of complex numbers, inquiring what properties they have, then asking whether and how they ``exist'' is somewhat beside the point.\pa

There are ``formalists'' who try to reject the ``mathematical objects'' entirely. They say that mathematics is only a performance, done according to certain rules: one knows, for example, what an assertion is and when one has completed a proof of such an assertion making it into a theorem, but that is all: this theorem does not teach us anything about any kind of ``objects''.\pa

This attitude seems to me, however, to be in contrast with the nature of the mathematician as an explorer. What
the mathematician is usually after is not doing some ``performance'' but learning about something. While the jogger must actually run to get any benefit, a mathematician would usually be content if (s)he is sure a proof can
be carried out. In fact, what mathematicians working in the above ``formalistic'' spirit are doing is overwhelmingly not ``performing'' but exploring. They explore the mathematical system imitating the part of the ``real world'' consisting of a crucial ingredient of the work of ordinary mathematicians, forming a an important subject-matter of Mathematical Logic.\pa

One may insist that just as a realm of ``objects'' needs a ``logic'' to be spoken about, so a ``logic'' is such only if it speaks of something. Otherwise it is no genuine logic at all, just some system $=$ universe of
discourse, to be spoken about by its own ``logic''.%
\footnote{One may note a somewhat similar ``duality'' about the World/Reality contrasted with our Knowledge.
In particular, we cannot say anything about the World but by our Knowledge, however dependent on the particulars of our human nature it is when viewed as part of the World.}
\pa

One often hears people debating about the essence (``way of existence'') of mathematical entities, as if somehow trying to ``find'' them in nature or our minds or society or the lore of philosophy etc. Sometimes the question raised is where these entities (or notions) ``exist'', and sometimes they are even rejected if not satisfactorily ``found''.\pa

In some sense, one thus follows the shadows of the Ancient Greeks who ``naturally'' viewed their mathematical objects as having a certain ``concrete'' nature close enough to experience or visual intuition (so Euclid did not consider negative numbers, and only grudgingly ``ratios'', i.e.\ positive real noninteger numbers). But Hermann Weyl (and Intuitionism) do that not with the ``naive'' notions of the Ancient Greeks, but with things like the members of $\bN$ and $\bR$ in the 20th century sense, rejecting what they cannot convince themselves to have ``found''.

\subsection{Claim}
One must choose sides: either one demands rigorous proofs, or one views customary rigor as a superfluous indulgence and is content with experimental or heuristic arguments, considered sufficient to obtain knowledge. Also, the role of computers in proofs is a kind of paradigm-change.\pa
\Ch
It may be said that the defining character of mathematics is that logical reasoning is the main path of investigation, either in an explicit orderly manner like the way it has been done in the mathematical literature since Euclid, or in a more implicit way customary in schools and in many past mathematical cultures. Certainly the validity of the usual methods to do arithmetics with numbers bigger than $10$, or of the formula for the area of a rectangle, have been usually justified by reasoning, although the Babylonians, e.g., did not find it necessary to keep records of such justifications.\pa

Still, Even with a system providing ``proofs'', there may be two ``ways of conduct''. One may make the judgment whether one really has a proof during the action itself, so to speak, based, maybe, on taste and the convincing power of the argument, or one may specify in advance what a proof should be, so that the action of mathematical reasoning can itself be made into a (``concrete'') mathematical system comparable to Magnitudes or Geometry. The latter way of conduct is preferred in modern mathematics -- one says then that the mathematics can be ``formalized''. This way also guarantees (as Bourbaki mention in their historical notes) a divorce ``in the acting'' between the ``mathematics'' and any kind of intuition or philosophical position. The former way of conduct was the only way in mathematics before the 20th century, notably also in Euclid (especially with respect to  order, e.g.\ judging when for points $A$, $B$, $C$ on a straight line $AB=AC+BC$ and when $AB=AC-BC$ or $AB=BC-AC$).\pa

This is also the way of physicists. They strive to calculate from theory to get, very often numerical but sometimes functional etc.\ results which may be compared with experiment or used for predictions. They take the liberty to reason in ways not (yet?) sanctioned by rigorous mathematics, and make the judgment whether the way of calculation is correct ``on the spot'', so to speak. Their ``only'' risk would be getting two different values for the same quantity, and even then they may decide ``on the spot'' to reject one of the
two ways of calculation.%
\footnote{This often extends to physicists' attitude to what constitutes the physics theory itself -- it too need not have a logical basis ``given in advance'', one rather decides ``on the spot'' what the theory says. In particular, this seems to be the way most physicists approach quantum theory.}
\footnote{Such mode of reasoning may very well reject a proof totally rigorous by the usual standards of ``rigorous mathematics'' if it seems ``too fancy''. Thus one may not acknowledge a contradiction that ``rigorous mathematics'' does. In particular, one may not find that a contradiction leads us to prove everything -- a basic theorem in usual Mathematical Logic.}
A piece of their mathematics would become ``rigorous'' the moment it could be cast into, say, Set Theory, which should happen if it (i.e.\ the criterion for the correctness of the reasoning) can somehow be ``formalized'' in the above sense.\pa

Let me clarify that. Sometimes in ``non-rigorous'' reasoning about mathematics (say, by physicists) one manipulates, say, integrals, in the face of the lack of any general theorem that says such manipulations are correct, and even the existence of counterexamples (a famous case is Feynman's path integrals, which fly in the face of known ``rigorous'' ways to define integrals). As long as one follows the first ``way of conduct'' above, one does trust the calculation (which naturally has avoided the counterexamples). Yet no one knows how to offer definitions and theorems which will give correct and non-contradictory results whenever anyone uses them. But if such definitions and theorems are found in the future, then this piece of reasoning would become ``rigorous''.\pa

Nowadays, with the availability of computers, one can often do experiments with mathematical objects and obtain conclusions ``with high probability'' similarly to ``real world'' science. For example, engineers have total conviction in a numerical method (coming with a plausible ``rationale'') to solve a type of partial differential equations after they had tested it on a sufficiently comprehensive sample of cases, and tests for primality of a natural number with a highly probable answer \cite{Rabin} are famous. Also, physics-like, more or less ``heuristic'' arguments succeed in ``establishing'' many mathematical assertions for which a more ``rigorous'' proof is elusive. The truth of the ``fact'' thus ``proved'' may be then judged either to be ``known'', or to be known with high certainty as in ``real world'' science.\pa

One should note, however, that as in ``real world'' science the sampling must be appropriate -- as said above, the set of natural numbers up to some $M$, however big, is a poor sample for the set of all natural numbers, thus testing, say, the Riemann Hypothesis up to some $M$ might not be viewed as a convincing sampling, even if one wishes to ``imitate'' natural sciences.\pa

Also, although there can be no questioning of the above-mentioned conviction of engineers, one may ``test'' whether, say, a probabilistic conclusion that a given $n$ is prime has the same force as a rigorous proof by imagining a match (as in sports) between them, namely,
wondering what would be said if, say, a rigorous proof and a probabilistic (or heuristic) conclusion led together to a contradiction.\pa

If one does not seek proofs that in principle give only ``highly probable knowledge'', but wishes to ``prove'', in the old-new sense of mathematics since Euclid -- prove so as to ``silence'' all objection -- then people's trust in ``not rigorous enough proofs'' of the great mathematicians of the past and of physicists and those who follow their pattern, depends on one's following the first ``way of conduct'' above, that of making the judgment whether one really has a proof during the action itself, based on taste and the convincing power of the argument, rather than specifying in advance what a proof should be, making the action of mathematical reasoning formalizable.\pa

In proofs that can be formalized%
\footnote{Although, of course, one does not care to explicitly write the formalization, this being a ``waste of time'', since one knows that there is no question this can be done if one wishes so and has the time and patience}
all deduction may be judged by exact rules, given in advance, valid with no exception. If one tries to cast the ``not rigorous enough'' arguments in such a form, one would find that task impossible -- the attempted rules will have exceptions, will be vague in the sense that one still has to decide what they dictate, or will lose their persuasive power. Thus, just choosing the second way of conduct, that of formalizable reasoning, makes the idea of ``non-rigorous proof'' evaporate. Then any disagreement boils down to a precise difference about which axioms to use (or, say, which rules of inference) and in particular, as Yuri Manin has pointed out, then one can discern unambiguously whether a purported ``proof'' is a rigorous proof (so then a mathematician who believes (s)he has a proof can be found wrong only due to a genuine ``mistake'' and not because ``standards of rigor have changed'').\pa

Also, it might be helpful to stress the nature of a proof as something ultimately formalizable thus capable of being viewed as a mathematical object itself. In this spirit, one should distinguish between the ``formal'' force of the proof and the ``real world'' question whether I really have a proof, so whether I really know that the theorem is true: have I not made a mistake? can I trust my memory? can I trust this book's claim that so and so was proved? can I trust the hardware and software of a computer that I used in the
proof? etc.\pa

Note, that G\"odel's Incompleteness Theorem is sometimes taken as denying mathematics a totally formal/mechanical/suitable for computers character (which it might have had Hilbert's program succeeded), guaranteeing its human/creative nature. Another cluster of mathematical ideas -- Computational Complexity Theory -- as if counters this. A theorem of length $\ell$ whose shortest proof is of length more than, say, $\ell^{10}$ is practically not a theorem for us, and proving theorems in the sense of finding a proof of length not more than $\ell^{10}$ (if there is such a proof) is ``just'' an $\textrm{NP}$-complete problem, equivalent to any other $\textrm{NP}$-complete problem (such as coloring a graph by $3$ colors) if computation that takes polynomial time is considered ``easy'', something to do quickly by machine, in a sense ``trivial''. In this way mathematics becomes an afternoon riddle -- coloring a graph, and only the ``hope'' that $\textrm{P}\ne\textrm{NP}$ prevents it from being trivial altogether in the just mentioned sense. But one might ``answer'' this (and let the pendulum swing back to ``human mathematics'') by noting that all this refers to proof as a formal (indeed mathematical) entity, while the proofs that mathematicians think, discuss and write are something different, ``human and unmathematical''. For example, an author of a mathematical article, or an editor of a journal, endowed by results of Computational Complexity Theory with an efficient formal way to check proofs written formally, still has to check whether the formal version is really a rendering of the ``informal'' ideas (so the formally-written proof may be found wrong though the ideas are correct).\pa

In ``rigorous'' mathematics, the task of a proof, which makes it a proof, is to rule out the assertion being false (and thus to make any attempt to prove its negation ``futile and ridiculous''). The all-importance of adequate rigor has been pointed out to mathematicians throughout their history by many plausible, even ``obvious'' assertions which were subsequently proved false (an example was the ``obvious'' conviction that every continuous function must be differentiable ``almost at all points''). In this vein, one may ``forgive'' the way of physicists and of those who make heuristic, physicists-like arguments in mathematics proper, as a taking of a calculated risk. Note that for physicists, merely deriving conclusions from any mathematical model, inevitably ``neglecting most of what happens in the universe and replacing what is inexactly measured with mathematical magnitudes'' incorporates such calculated risk.\pa

Of course, a ``worthy'' mathematical proof does more than just assuring us that the assertion being false has been ruled out. It shows us the ``wisdom'' behind that. It, as Yuri Manin says, makes us wiser. Also, these ``rigorous'' proofs carry an undeniable part  of the claim of mathematics to be beautiful and appealing. (In fact, very often most of the wisdom, and much of the beauty, lies in the proof and in understanding it, not in the theorem as such.)%
\footnote{To quote Freeman Dyson (in the context of subjects like Chaos Theory where one sees lots of numerical graphs and pictures), rigorous theorems are the best way to give a subject intellectual depth and precision.}\pa

In conclusion we are left with both practices -- that of engineers and physicists relying on testing suitable samples and on ``heuristics'', arguably quite adequate to keep them (and those who choose to treat mathematical objects their way) on the right track for their needs (and certainly not lacking as a source of wisdom, beauty and appeal) and that of ``rigorous'' mathematicians seeking rigorous (and preferably enlightening) proofs -- each having its proper legitimate place.\pa

\subsection{Claim}
Sometimes a set of ``real world'' entities has a structure of a set of magnitudes which can be added -- in fact, by choosing a zero point and a unit of measurement we can equate them with the set of numbers (reals, for ``continuous'' entities, or integers, for ``discrete''). In other cases they are naturally endowed only with an ordering -- one can say which is bigger -- mathematically they have the structure of a topological line, say. Yet, in the latter case we strive to ``quantify'' them -- turn them into magnitudes as in the former case, thus enriching their mathematical relevance.\pa
\Ch
Of course, an automorphism of the order structure, i.e.\ a one-one transformation of, say, the topological line onto itself that preserves its order structure, hence preserves everything significant in the system if only the order made sense there, will, in general, transform the imposed addition structure into a totally different one. Hence, it will transform any consequences derived from the ``quantified'' system into ``arbitrarily'' different ones. Such consequences should thus be considered, in general, as arbitrary and irrelevant. Thus, if one just ``quantifies'' a system whose natural mathematical model is an ordered -- topological -- line (let alone where there is even no natural order), and treats these numbers as if adding, averaging them etc.\ has meaning, one often abuses, rather then uses, mathematics.


\begin{thebibliography}{MM1}
\bibitem[We]{Weyl}
Hermann~Weyl, {writings, in a constructivist/intuitionist mood, about which things in mathematics really exist and are legitimate}.
\end{thebibliography}
\end{document}